\documentclass{amsart}
\usepackage{amsmath,amssymb,amscd,amsthm,amsfonts,graphicx,hyperref}

\newcommand{\R}{\mathbb R}
\newcommand{\C}{\mathcal C}
\newcommand{\F}{\mathcal F}

\newtheorem{theorem}{Theorem}
\newtheorem{lemma}{Lemma}

\title{Paths on the doubly covered region of a covering of the plane by unit discs}
\author{Edgardo Rold\'an-Pensado}
\date{}

\begin{document}

\begin{abstract}
Given a covering of the plane by closed unit discs $\F$ and two points $A$ and $B$ in the region doubly covered by $\F$, what is the length of the shortest path connecting them that stays within the doubly covered region? This is a problem of G. Fejes-T\'oth and he conjectured that if the distance between $A$ and $B$ is $d$, then the length of this path is at most $\sqrt 2 d+O(1)$.
In this paper we give a bound of $2.78 d+O(1)$.
\end{abstract}

\maketitle

\section{Introduction}

Let $\F$ be a locally finite covering of the plane $\R^2$ by closed unit discs.
The doubly covered region of $\R^2$ by $\F$ consists of the sets of points $x\in\R^2$ that are contained in at least two elements of $\F$.

A few years ago Gabor Fejes-T\'oth posed the following question: If $A,B\in\R^2$ are two points at distance $d$ apart and contained in the doubly covered region of $\R^2$, what is the length of the shortest path $\gamma$ joining $A$ and $B$ that is completely contained in the doubly covered region?

In a sense, this is the dual of a problem by Laszlo Fejes T\'oth about the length of a path avoiding a packing of discs \cite{Fej1993}.

Gabor Fejes-T\'oth conjectured that when the centres of the circles in $\F$ form a unit square lattice, the length of $\gamma$ is maximal.
For any two points $A$ and $B$ in this example, $\lvert\gamma\rvert=\sqrt 2 d + O(1)$.
A general upper bound of $\pi d+O(1)$ for $\lvert\gamma\rvert$ is not difficult to obtain.

Baggett and Bezdek proved in \cite{BB2003} that when the centres of the circles of $\F$ form a lattice, then the unit square lattice is indeed the extreme case.

In this short note we give an upper bound for any locally finite covering.

\begin{theorem}
$$\lvert\gamma\rvert\leq\left(\frac{\pi}{3}+\sqrt 3\right)d+O(1)\approx 2.77925 d+O(1).$$
\end{theorem}

This is still closer to $\pi d$ than to $\sqrt d$, but we hope that our methods may help others to continue improving this bound.

\section{Proof of Theorem 1}

Let $L$ be the segment joining $A$ and $B$ and $\mathcal G\subset\F$ be a minimal sub-cover of $L$.
Suppose that $L$ is horizontal and that $A$ is to the left of $B$.
The elements of $\mathcal G$ can be ordered as $\C_1,\dots,\C_n$ such that if $i<j$ then the centre $O_i$ of $\C_i$ is to the left of the centre $O_j$ of $\C_j$.

We may assume that $A$ is the leftmost point of the intersection of the line $AB$ with $\C_1$ and that $B$ is the rightmost point of the intersection of $AB$ with $\C_n$.
If this were not the case then we may extend the length of $L$ by at most $4$ and move each one of $A$ and $B$ through paths of length at most $2+\pi$ so that they end up in this way.
This contributes a term $O(1)$ to the length of the curve we find with respect to the original length of the segment $L$.

Since the family $\F$ is locally finite, every point in the boundary of a circle $\C\in\mathcal G$ is doubly covered.
For $\C\in\mathcal G$, define $\C'$ as the closure of a connected component of $\C\setminus L$ that does not contain the centre of $\C$.
Since $\mathcal G$ is a minimal covering of $L$ then $\C_i'\cap \C_j'\neq\emptyset$ if and only if $\lvert i-j\rvert\leq 1$.

Let $A_i$ and $B_i$ be the points of intersection of the boundary of $\C_i$ with $L$ such that $A_i$ is to the left of $B_i$ and let $M_i$ be the midpoint of the segment $A_{i+1}B_i$ for $0<i<n$, $M_0=A_1$ and $M_n=B_n$.
Note that $\lvert M_iM_{i-1}\rvert\ge\lvert A_iB_i\rvert/2$.

\begin{figure}
\centering
\includegraphics[scale=0.75]{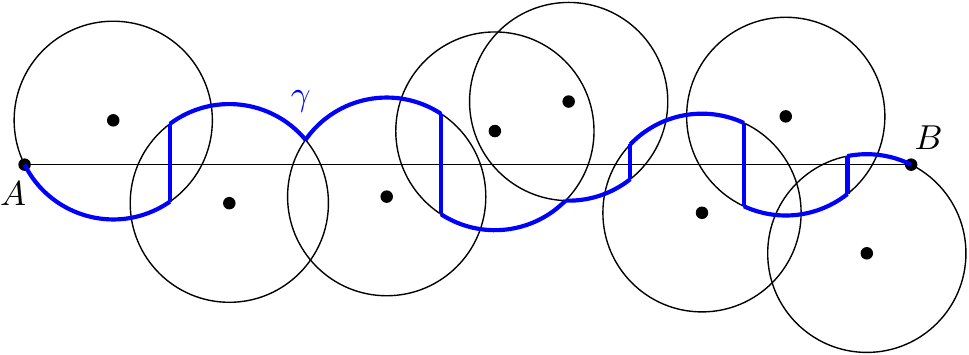}
\caption{The path $\gamma$.}
\label{fig:gamma}
\end{figure}

Now we construct the path $\gamma$ with the algorithm below, an example is shown in Figure \ref{fig:gamma}.

The path starts at $A=M_0$.
Assuming the path has been constructed up to $M_i$, let $k$ be the largest integer such that the sets $\C_{i+1}',\dots,\C_{i+k}'$ are all on the same side of $L$.
Without loss of generality we assume they are all above $L$.
The path then continues vertically upwards until it reaches the boundary of $\C_{i-1}$.
From here it continues towards the right while staying contained in the boundary of $\bigcup_{j=1}^k\C_{i+j}'$ until it is  vertically aligned with $M_{i+k}$.
Finally it goes vertically downwards until it reaches $M_{i+k}$.
This is repeated until the path ends at $B=M_n$.

It is easy to see that $\gamma$ is well defined and completely contained in the doubly covered region, however it is not so easy to directly bound its length.
This path could be shortened by taking diagonal lines instead of vertical ones, but these two paths coincide in the extreme case.

\begin{figure}
\centering
\includegraphics[scale=0.75]{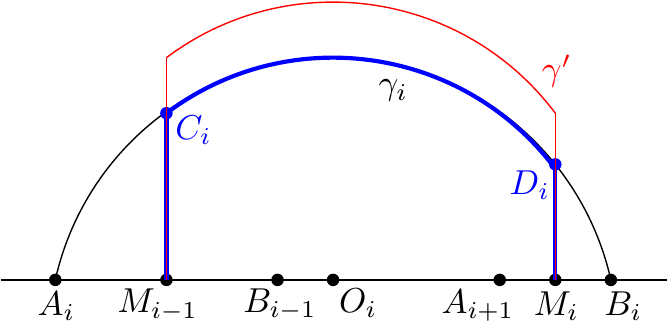}
\caption{The paths $\gamma_i$ and $\gamma'$.}
\label{fig:gammai}
\end{figure}

In order to bound the length of this path, for every $\C_i\in\mathcal G$ we construct a new path $\gamma_i$.
The paths $\gamma_i$ may not be contained in the doubly covered region but they satisfy $\lvert\gamma\rvert\leq\sum\lvert\gamma_i\rvert$.
Assume for simplicity that $\C_i'$ is above $L$, then the path $\gamma_i$ starts at $M_{i-1}$, goes upwards until it intersects the boundary of $\C_i'$ at a point $C_i$, then goes to the right staying contained in the boundary of $\C_i'$ until it is vertically aligned with $M_i$ at the point $D_i$, and finally goes downwards until it reaches $M_i$. This path is shown in Figure \ref{fig:gammai}.
It is not difficult to see that $\lvert\gamma\rvert\leq\sum\lvert\gamma_i\rvert$.

Now it is enough to prove the following.
\begin{lemma}\label{lem:lemma}
$$\lvert\gamma_i\rvert\leq\left(\frac{\pi}{3}+\sqrt 3\right)\lvert M_{i-1}M_i\rvert.$$
\end{lemma}

We may assume that $A_iB_i$ is the diameter of $\C_i$, otherwise let $\C'$ be a circle with diameter $A_iB_i$ and $\gamma'$ be the curve defined on $\C'$ in the same way that $\gamma_i$ is defined on $\C_i$.
This new curve clearly has larger length than $\gamma_i$ (see Figure \ref{fig:gammai}).

Let $\alpha=\angle C_iO_iA_i$ and $\beta=\angle B_iO_iD_i$, then
\begin{align*}
\lvert\gamma_i\rvert=&\sin(\alpha)+\sin(\beta)+\pi-\alpha-\beta,\\
\lvert M_{i-1}M_i\rvert=&\cos(\alpha)+\cos(\beta).
\end{align*}
Recall also that $\lvert M_{i-1}M_i\rvert\ge \lvert A_iB_i\rvert$.

If we fix the value of $\lvert M_{i-1}M_i\rvert$, then by using Lagrange multipliers we obtain that $\lvert\gamma_i\rvert$ is maximum when $\alpha=0$, $\beta=0$ or $\alpha=\beta$.
Since the cases $\alpha=0$ and $\beta=0$ are symmetrical, we have essentially two cases.
Now we only need to determine the maximum of
$$\max\left\{\frac{\sin(\alpha)+\pi-\alpha}{\cos(\alpha)+1},\frac{2\sin(\alpha)+\pi-2\alpha}{2\cos(\alpha)}\right\}$$
as a function of $\alpha$.
This occurs when $\alpha=\frac{\pi}{3}$ and corresponds to the case $\alpha=\beta$. This gives
$$\frac{\lvert\gamma_i\rvert}{\lvert M_{i-1}M_i\rvert}\leq\frac{\pi}{3}+\sqrt 3,$$
which proves the lemma.

\section{Limitations of our approach}

The method we use only considers discs that intersect the segment $AB$, we construct a path contained in the boundary of these discs and in their doubly covered region.
Considering only these circles it is impossible to obtain the bound Fejes-T\'oth conjectured.

Below we construct an example considering only these circles such that $\lvert\gamma\rvert>\sqrt 2 d$.

Let $\C_1$ and $\C_2$ be intersecting circles and assume $A\in\partial\C_1\setminus C_2$ and $B\in\partial\C_2\setminus C_1$.
If we only allow the path $\gamma$ to be in the intersection of the circles and their boundary, then there are only $4$ possible choices for $\gamma$.
This is a simplified version of our problem, as the boundary of a circle in $\F$ is always doubly covered.
But because we only consider two circles, this new problem can be solved precisely to obtain $\lvert\gamma\rvert\le 1.58003 d$.
Here the number $1.58003$ is an approximation obtained as the maximum of a function involving trigonometric functions.

\begin{figure}
\centering
\includegraphics[scale=0.75]{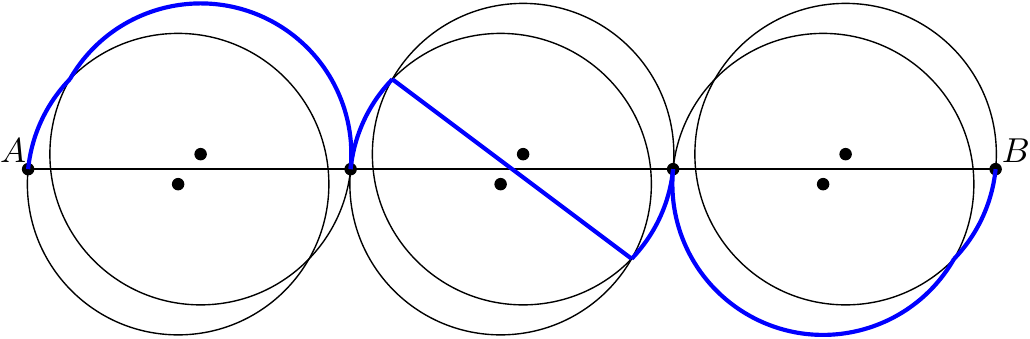}
\caption{A path of length $1.58003 d$.}
\label{fig:example}
\end{figure}

The example that gives this value can be extended to families with more circles and have arbitrarily large $d$.
This is done by using several translated copies of this example (see fig. \ref{fig:example}).

Thus, it is essential to consider elements of $\mathcal F$ further away from the segment $AB$.
However an improvement on our algorithm would surely lower our bound.

\section{Remarks}

To obtain our bound we consider only $2$ discs at a time, this is precisely what is done in Lemma \ref{lem:lemma}.
A this bound could potentially be lowered if we were to consider $3$ or more discs instead.
The algorithm would have to be modified or changed to take advantage of this.

We tried this approach together with Alexey Garber, however we were unable to obtain any substantial improvements.
A problem we encountered was that the functions to be minimised became extremely complicated.

\end{document}